\documentclass[journal, 10pt]{IEEEtran}
\pdfoutput=1
\usepackage{amsmath}
\usepackage{graphicx}
\usepackage{amsfonts,amssymb}
\graphicspath{ {images/} {Figures/}}
\usepackage{blindtext, amsmath}
\usepackage[backend=bibtex,sorting=none,]{biblatex}
\usepackage{color}
\hfuzz=\maxdimen
\ifCLASSINFOpdf
\else
\fi
\begin{document}
\title{Robust Beamforming Against Direction-of-Arrival Mismatch Using Subspace-Constrained Diagonal Loading}
\author{Yueh-Ting Tsai, Borching Su, Yu Tsao, and Syu-Siang Wang}

\maketitle

\begin{abstract}
In this study, a new subspace-constrained diagonal loading (SSC-DL) method is presented for robust beamforming against the issue of a mismatched direction of arrival (DoA), based on an extension to the well known diagonal loading (DL) technique.
One important difference of the proposed SSC-DL from conventional DL is that it imposes an additional constraint to restrict the optimal weight vector within a subspace whose basis vectors are determined by a number of angles neighboring to the estimated DoA.
Unlike many existing methods which resort to a beamwidth expansion, the weight vector produced by SSC-DL has a relatively small beamwidth around the DoA of the target signal.
Yet, the SSC-DL beamformer has a great interference suppression level, thereby achieving an improved overall SINR performance.
Simulation results suggest the proposed method has a near-to-optimal SINR performance.

\end{abstract}


\IEEEpeerreviewmaketitle

\section{Introduction}
Beamforming is a signal processing technique that utilizes the spatial information from a set of receivers to attain the signal with satisfactory quality. It has been used in different applications, wireless communications, acoustics, etc.    As a data-dependent beamformer, the classical Capon method \cite{capon1969} is well known to achieve this goal and maximizes the output signal-to-interference-plus-noise (SINR) level using the minimum variance distortionless response (MVDR) criterion. 
However, in a presence of a steering vector error caused by, for example, a mismatch in a direction of arrival (DoA), the performance of the beamformer seriously degraded since the desired signal is mistreated as an interference and suffers a severe attenuation. 
In order to resolve this issue, numerous approaches have been proposed for dealing with the steering vector errors. 
A class of beamformers known as linear constrained minimum variance (LCMV) were proposed to mitigate this issue by imposing additional linear constraints that increase the beamwidth of the mainlobe \cite{Vural1977,Takao1976}.
Another popular class of beamformers, called diagonal loading (DL), addressed the issue by adding a scalar matrix to the covariance matrix \cite{bd1998,li2003}.
DL methods are known to provide robustness not only to the DoA mismatch issue but also to general types of steering vector errors. 
One well-known drawback of the DL methods is that a universally accepted way of choosing the diagonal loading factor is still lacking. 
Nevertheless, in \cite{vincent2004}, an optimal DL factor that maximizes SINR was derived with some approximations. 
Further, in \cite{ma2003,du2010}, automatic determination of DL factor values was presented.
In order to enhance the performance, many extensions to DL \cite{sq1999,vorobyov2003, Lorenz2005,ShahbazPanahi2003,chen2007} were proposed with an idea to force the magnitude response of all steering vectors in a sphere set centered at the desired steering vector to exceed or equal to unity while minimizing the output variance. 
In \cite{Feldman1994}, a method is proposed by projecting the steering vector onto the signal-plus-interference subspace. More recently, methods that based on complex quaternion processes \cite{tao2014}, based on the artificial immune system \cite{kiong2014}, based on interference-plus-noise covariance matrix reconstruction \cite{huang2015}, based on an interference subspace \cite{mao2015} were also proposed.

In this study, we consider only the steering vector error caused by DoA mismatches and propose a new extension to DL method called subspace-constrained diagonal loading (SSC-DL). 
The difference between SSC-DL and conventional DL is that the search space of weight vectors is constrained in a subspace determined by the signal covariance matrix and the possible range of the desired signal's DoA. 
An important feature of the proposed method over other DL methods is that the additional constraints are added to ensure an interference suppression capability of the beamformer instead of working only on broadening the beamwidth.
Such constraints turn out to be helpful in improving SINR performance. 
The simulation results show the effectiveness of the proposed SSC-DL approach in suppressing interference signals.

The rest of this paper is organized as follows. 
Section II introduces the MVDR beamforming system and the DoA mismatch issue. 
Section III presents the proposed SSC-DL algorithm. 
Section VI presents the derivation of the diagonal loading factor used in the proposed SSC-DL algorithm.
Section \ref{sec:sim} presents the simulation setup and results. 
Finally, conclusions are provided in Section \ref{sec:conclusion}.
Boldface lowercase letters represent column vectors; boldface uppercase letters represent matrices. {Superscripts $^T$, $^H$, $^\perp$, and $E[\cdot]$ represent matrix transpose, Hermitian-transpose operations, and subspace orthogonal complement, and expectation operation, respectively.}

\section{MVDR Beamformers}
In this section, we first introduce the beamforming system, MVDR beamformer, and the DL method. Next, we describe the DoA mismatch issue. 

\subsection{The Beamforming System}
\begin{figure}[htp]
\centering \centerline{
\includegraphics[width=0.5\textwidth,clip]{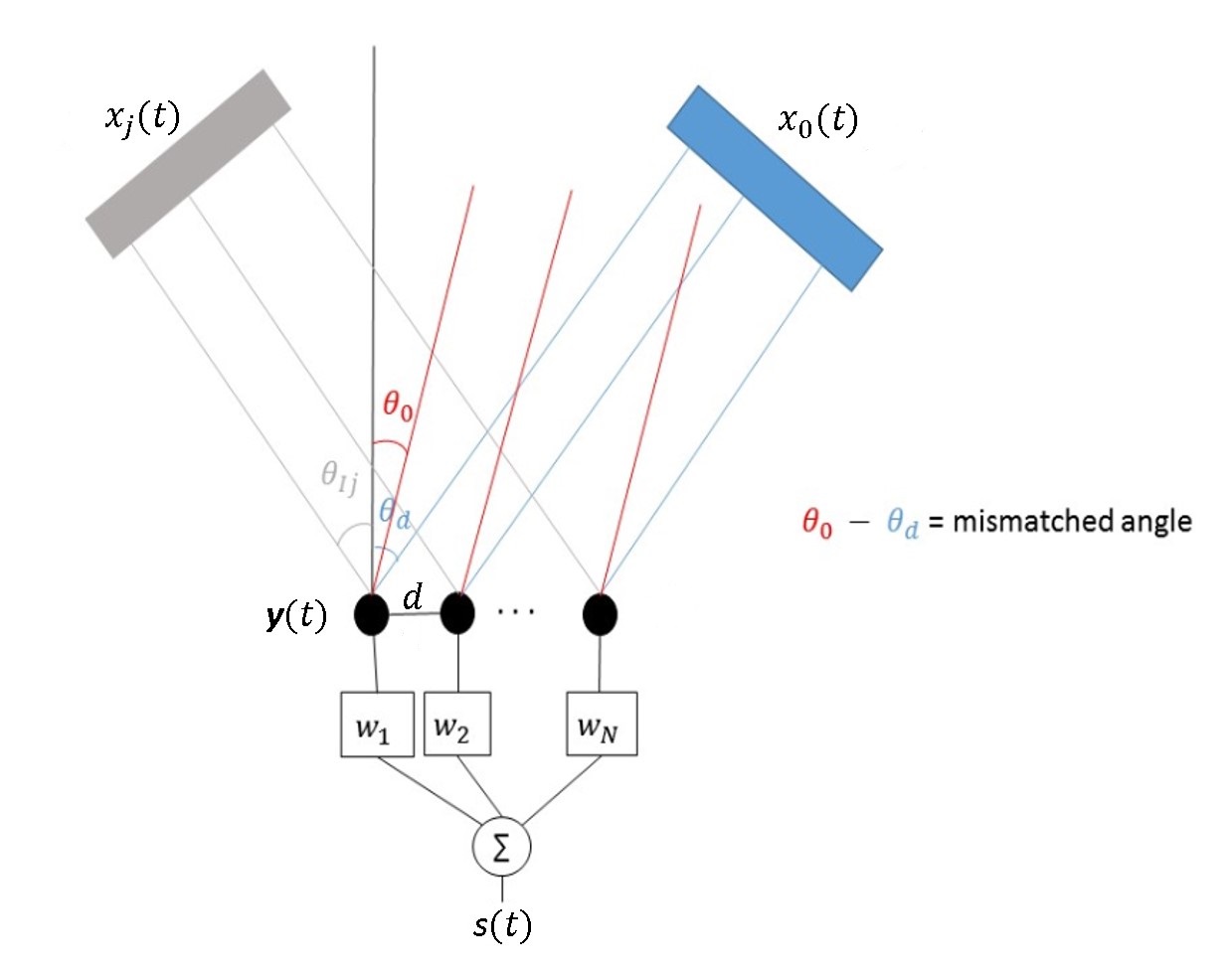}}
\caption{system model for robust beamforming}\label{fig:system}
\end{figure}

With figure \ref{fig:system}, consider a beamforming system comprising $N$ sensors in a uniform linear array (ULA) with an inter-element spacing $d$. 
The estimated DoA of the desired signal is $\theta_0$, and the actual DoA of the desired signal is $\theta_d$. 
When $\theta_0$ is not equal to $\theta_d$, the difference $\theta_0-\theta_d$ is called the DoA mismatch. 
We denote $\mathbf{a}(\theta)\in\mathbb{C}^N$ as the steering vector toward direction $\theta$ \cite{bd1998}: \begin{align}\label{eq2}
\mathbf{a}(\theta)=[1,e^{j(2\pi/\lambda)d\sin\theta},\cdots,e^{j(N-1)(2\pi/\lambda)d\sin\theta}]^{T} 
\end{align}
where $\lambda$ is the wavelength. 
Assume there are $J$ interference signals $x_j(t), j=1,..., J$, and the $j$th interference signal is arriving from angle $\theta_{Ij}$, {where $t\in{\cal R}$ denotes time.}
Then, the received vector signal  {at time $t$}, $\mathbf{y}(t)=[y_1(t), y_2(t),\cdots,y_N(t)]^T$ can be expressed as
\begin{equation}
\begin{aligned}\label{eq:yt}
\mathbf{y}(t)&=x_0(t)\mathbf{a}(\theta_d)+\sum_{j=1}^J x_{j}(t)\mathbf{a}(\theta_{Ij})+\mathbf{n}(t)\\
&=x_0(t)\mathbf{a}(\theta_d)+{\mathbf i}(t)+{\mathbf n}(t)
\end{aligned}
\end{equation}
where $x_0(t)$ is the desired signal and {${\bf n}(t)$ denotes additive white Gaussian noise whose entries are i.i.d. complex Gaussian random variables $\sim {\cal CN}(0, \sigma_n^2)$. The signals $x_0(t)$ and $x_{j}(t), j=1,..., J$ are assumed to be statistically independent zero-mean wide-sense stationary processes, with $E[|x_0(t)|^2] = \sigma^2_{x}$ and $E[|x_{j}(t)|^2] = \sigma^2_{I_{j}}$.} 
The desired signal $x_0(t)$ and the interference signals $x_j(t)$ are assumed to be uncorrelated.

The beamforming technique uses a weight vector $\mathbf{w}=[w_0, w_1,\cdots,w_{N-1}]^T$ to process the received signal ${\bf y}(t)$ in an attempt to suppress the interference component present in ${\bf y}(t)$ while preserving the desired signal $x_0(t)$.
The performance of the processed signal $  s(t) = {\bf w}^H {\bf y}(t)$ can be evaluated by the signal-to-interference-plus-noise ratio (SINR), defined as
\begin{align}\label{eq:SINR}
SINR=\frac{P\left | \mathbf{w}^H {\bf a}(\theta_d) \right |^2}{\mathbf{w}^H\mathbf{R}_{I+n}\mathbf{w}}
\end{align}
where $P$ = $E[|x_0(t)|^2]$ is the desired signal power and $\mathbf{R}_{I+n}$ = $\mathbf{R}_{I}$ + $\mathbf{R}_{n}$ is the sum of the interference-covariance matrix {$\mathbf{R}_{I}$ = $E[\mathbf{i}(t)\mathbf{i}(t)^H]$, with ${\bf i}(t)$ defined as ${\bf i}(t)=\sum_{j=1}^J x_{j}(t)\mathbf{a}(\theta_{Ij})$}, and noise-covariance matrix $\mathbf{R}_{n}$ = $E[\mathbf{n}(t)\mathbf{n}(t)^H] = \sigma_n^2{\bf I}_n$.

\subsection{MVDR Beamformer}
The MVDR beamformer estimates the weight vector $\mathbf{w}$ by solving the following optimization problem:
\begin{align}\label{eq:MVDR-RI+n}
&\min_{\mathbf{w}}\mathbf{w}^{H}\mathbf{R}_{I+n}\mathbf{w}  ~~~ s.t. ~ |\mathbf{w}^{H}\mathbf{a}(\theta_d)|=1 .
\end{align}
Generally, $\mathbf{R}_{I+n}$ is difficult to be precisely estimated, and thus the classical MVDR beamformer uses signal covariance matrix $\mathbf{R}_{y}=E[\mathbf{y}(t)\mathbf{y}(t)^H]$ to replace $\mathbf{R}_{I+n}$. {When $\theta_d$ is perfectly known and $|\mathbf{w}^{H}\mathbf{a}(\theta_d)|=1$ is satisfied}, it is shown \cite{capon1969} that (\ref{eq:MVDR-RI+n}) is equivalent to:
\begin{align}\label{eq:MVDR}
	&\min_{\mathbf{w}}\mathbf{w}^{H}{\mathbf{R}}_{y}\mathbf{w}  ~~~ s.t. ~ |\mathbf{w}^{H}\mathbf{a}(\theta_d)|=1.
\end{align}
{The equivalence of (\ref{eq:MVDR-RI+n}) and (\ref{eq:MVDR}) becomes clear by observing that
\begin{align}
	\mathbf{w}^{H}{\mathbf{R}}_{y}\mathbf{w} &= \mathbf{w}^{H}(P{\bf a}(\theta_d){\bf a}^H(\theta_d)+{\mathbf{R}}_{I+n})\mathbf{w} \nonumber\\
	&=P + \mathbf{w}^{H}{\mathbf{R}}_{I+n}\mathbf{w} .
\end{align}}
The optimal solution for (\ref{eq:MVDR}) can be expressed in a closed form:
\begin{align}\label{eq:MVDR-closedForm}
\mathbf{w}=\frac{{\mathbf{R}}_{y}^{-1}\mathbf{a}(\theta_d)}{\mathbf{a}^H(\theta_d){\mathbf{R}}_{y}^{-1}\mathbf{a}(\theta_d)}.
\end{align}
In practice, signal covariance matrix ${\bf R}_{ y}$ is obtained by averaging a number of snapshots $K$: 
\begin{align}\label{eq:hatRy}
{\hat{\mathbf{R}}}_{y} = \frac{1}{K}\sum_{k=1}^{K}\mathbf{y}(k)\mathbf{y}^{H}(k).
\end{align}

\subsection{The Diagonal Loading Method}
The MVDR beamformer is known to have poor performance if $\theta_d$ is not known accurately. 
When a DoA mismatch occurs, the desired signal may be mistreated as an interference and will be suppressed, thus decreasing the output SINR. 
The DL method \cite{bd1998, li2003} is known to be capable of mitigating the performance degradation caused by the DoA mismatch issue by including a scalar matrix in Eq. ($\ref{eq:MVDR}$):
\begin{align}\label{eq:DL}
{\bf w}_{DL}=\arg\min_{{\bf w}\in\mathbb{C}^N} {\bf w}^H ({\bf R}_y + \gamma{\bf I}){\bf w} ~~s.t.~ |{\bf w}^H {\bf a}(\theta_0)|=1
\end{align}
where $\gamma$ is a real number. 

The DL method is successful in avoiding the beamformer to misrecognize the desired signal as an interference; 
it is, however, limited in the capability of suppressing interference signals, leaving room for improvements in SINR performance.

\subsection{Problem Statement}
We consider a beamforming system with a DoA mismatch: the estimated DoA of the desired signal is $\theta_0$, which is unequal to the true DoA, $\theta_d$. Moreover, we assume the receiver knows a lower bound $\theta_1$ and an upper bound $\theta_2$ of the desired DoA, i.e., $\theta_1 < \theta_d < \theta_2$. {The DoAs of interference signals, $\theta_{Ij} $, are assumed not to be within the range $[\theta_1, \theta_2]$}. {We assume that the angles $\theta_d, \theta_{Ij},j=1,...,J$ are deterministic yet in general unknown to the receiver.} The goal is to use $\hat{\bf R}_y$ defined in (\ref{eq:hatRy}) and find a weight vector ${\bf w}$ that achieves a good SINR performance defined in (\ref{eq:SINR}). Without loss of generality, we assume all angles $\theta_d$, $\theta_0$, $\theta_1, \theta_2$, and $\theta_{Ij}$ are within the range $(-\pi/2,\pi/2)$.

\section{The Proposed SSC-DL Beamformer}
Although the MVDR with DL is capable of mitigating the unwanted suppression of the desired signal, there is still room to further improve the SINR performance in the presence of the DoA mismatch. 
In this section, we first describe the implementation steps of the proposed SSC-DL beamformer followed by an explanation of its advantages on interference suppression.

\subsection{Algorithm Description}
\begin{figure}[htp]
\centering \centerline{
\includegraphics[width=0.5\textwidth,height=5cm]{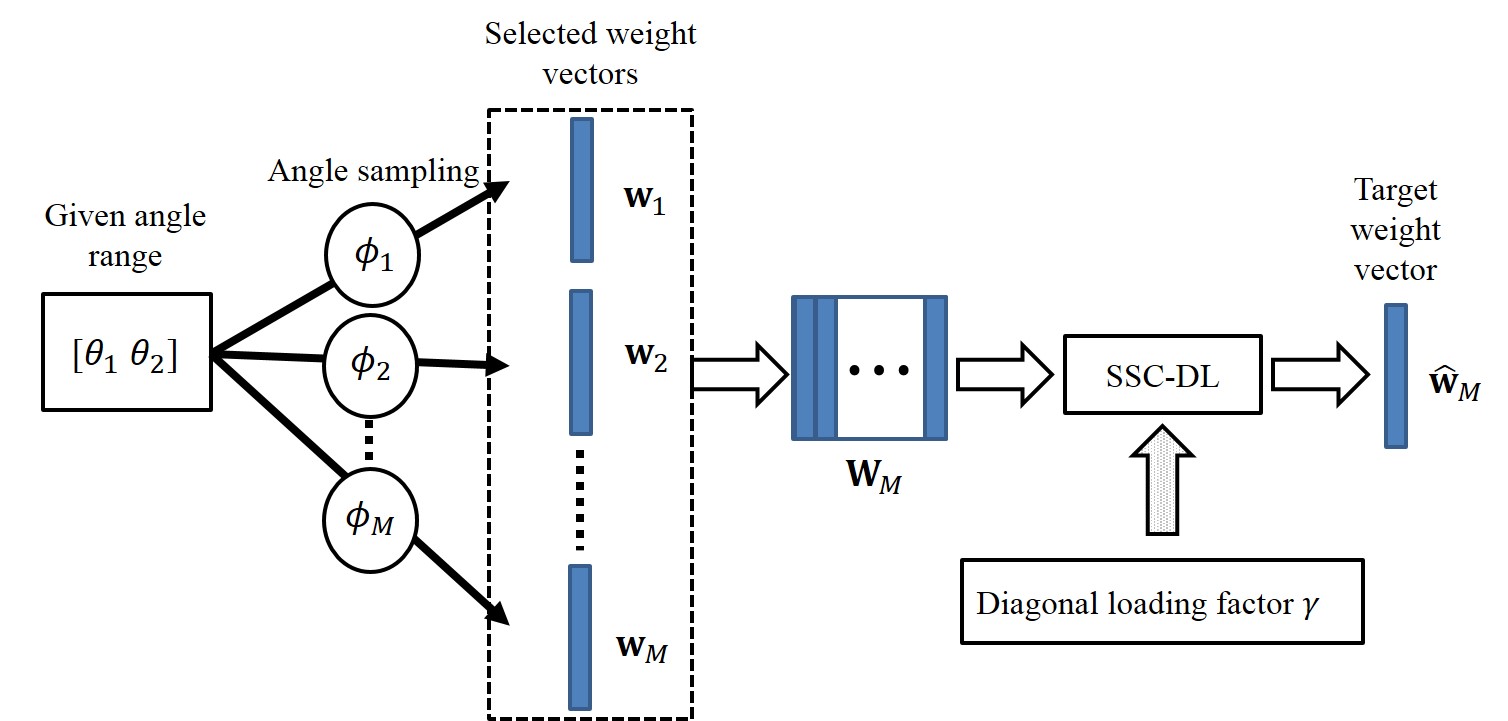}}
\caption{Flow Chart of SSC-DL}
\label{flsscdl}
\end{figure}
{Figure \ref{flsscdl} shows the flow chart of the proposed SSC-DL algorithm, which can be divided into three steps.} In the first step, SSC-DL chooses an integer $M<N$ and defines a set of $M$ angles neighboring to the estimated DoA, $\theta_0$, within the range $[\theta_1, \theta_2]$:
\begin{align}\label{eq:phi_m}
\phi_m = \sin^{-1}\left(\sin{\theta_1} + m \cdot \frac{\sin{\theta_2} - \sin{\theta_1}}{M-1}\right), 
\end{align}
$m=0,1,...M-1$. Note that $\theta_1=\phi_0<\phi_1<\cdots<\phi_{M-1}=\theta_2$. Then, for the $m$th angle $\phi_m$, we compute its corresponding weight vector according to (\ref{eq:MVDR-closedForm}):  
\begin{align}\label{eq:w_m}
{\bf w}_m = {  \frac{{\bf R}_y^{-1}{\bf a}(\phi_m)}{{\bf a}^H(\phi_m){\bf R}_y^{-1}{\bf a}(\phi_m)}},
\end{align}
which is the optimal weight vector in terms of SINR when the true DoA happens to be $\theta_d=\phi_m$.
In the second step, SSC-DL collects these $M$ vectors
and constructs an $M$-dimensional subspace of $\mathbb{C}^N$
\begin{equation}
{\cal W} = \left\{ {\bf W}_{M}{\bf z} | {\bf z}\in\mathbb{C}^M \right\}
\label{eq:subspaceW}
\end{equation}
where ${\bf W}_{M} = \left[\begin{array}{cccc}{\bf w}_1 & {\bf w}_2 & \cdots & {\bf w}_M \end{array}\right]$ is of size $N\times M.$
The third step is to compute the weight vector as:
\begin{align}\label{eq:VS}
{\bf \hat w}_{M} = \arg\min_{{\bf w}\in {\cal W}} {\bf w}^H ({\bf R}_y + \gamma{\bf I}){\bf w} ~~~s.t. ~|{\bf w}^H {\bf a}(\theta_0)|=1.
\end{align}
Here, the value $\gamma$ is chosen as 
\begin{equation} \gamma = -(\sigma_{n}^2 + P\cdot N),\label{eq:gamma}
\end{equation}
based on a derivation that will be later elaborated in Section \ref{sec:derivation}.
Using a technique of Lagrange multipliers, the optimization problem defined in (\ref{eq:VS}) can be shown to have a closed-form solution
\begin{align}\label{ssc-dl}
{\mathbf {\hat w}}_{M}=\frac{{\bf W}_{M}({\bf W}_{M}^{H}(\mathbf{R}_{y}+\gamma\mathbf{I}){\bf W}_{M})^{-1}{\bf W}_{M}^{H}\mathbf{a}(\theta_0)}{\mathbf{a}^{H}(\theta_0){\bf W}_{M}({\bf W}_{M}^{H}(\mathbf{R}_{y}+\gamma\mathbf{I}){\bf W}_{M})^{-1}{\bf W}_{M}^{H}\mathbf{a}(\theta_0)}.
\end{align}

\subsection{Choices of the subspace dimension $M$}\label{subsec:M}
The choice of the dimension of the subspace ${\cal W}$, $M$, could be any integer between $1$ and $N$. 
When $M=1$, the SSC-DL reduces to the  MVDR beamformer. 
On the other hand, when $M$ is chosen as $N$, the SSC-DL degenerates to the case of conventional DL method since the subspace ${\cal W}$ equals to the whole space $\mathbb{C}^M$ and there are effectively no constraints to the beamformer coefficients ${\bf w}$.

While we are not yet able to give a comprehensive criterion to an optimal choice of $M$ in this work, we provide here some rules of thumbs of the selection of the parameter. 
First of all, $M$ should be chosen to be smaller than or equal to $N-J$ since the subspace ${\cal W}$ should be designed to be ``almost orthogonal'' to the interference subspace which has a dimension of $J$.
On the other hand, $M$ should be chosen sufficiently large to cover all steering vectors covering all possible DoAs of the desired signal defined as the range $[\theta_1, \theta_2]$.
In Section \ref{subsec:SINRvsM}, the effect on the system performance of the parameter $M$ will be studied through numerical results.

\subsection{Interference Suppression}
\label{subsec:interf_suppr}
We give a qualitative explanation to why the proposed SSC-DL method possesses a capability of suppressing interference signals superior to the conventional DL method. 
The subspace ${\cal W}$ as defined in (\ref{eq:subspaceW}) can be interpreted as space that is ``almost orthogonal'' to all interfering steering vectors ${\bf a}(\theta_{Ij}), j=1,..., J$. It is also a space that the steering vector of the desired signal, ${\bf a}(\theta_d)$, with any possible angle $ \theta_d\in [\theta_1,\theta_2]$, would have most of its energy within.
Let ${\bf Q}_s$ and ${\bf Q}_n$ be matrices whose columns, respectively, form an orthonormal basis for the subspaces ${\cal W}$ and ${\cal W}^\perp$.
Then, the steering vector of an interfering signal will have most of its energy in ${\cal W}^\perp$:
\begin{equation} || {\bf Q}_s^H {\bf a}(\theta_{Ij})|| \ll ||{\bf Q}_n^H  {\bf a}(\theta_{Ij})|| \approx ||{\bf a}(\theta_{Ij})||,~ j=1, 2, ..., J. 
\nonumber \label{eq:QsllQn}
\end{equation}
Denote $\lambda_j$ as the ratio $\lambda_j = || {\bf Q}_s^H {\bf a}(\theta_{Ij})|| / ||{\bf Q}_n^H  {\bf a}(\theta_{Ij})||$. 
Then we have $\lambda_j \ll 1$. 
Further, noting that ${\bf \hat w} = {\bf Q}_s{\bf z}'$ for some ${\bf z}'$ satisfying $||{\bf \hat w}|| = ||{\bf z}'||$, we have
\begin{eqnarray}
 |{\bf a}^H(\theta_{Ij}){\bf w}| &=& |{\bf a}^H(\theta_{Ij}) {\bf Q}_s {\bf z'}|  \leq  || {\bf Q}_s^H {\bf a}(\theta_{Ij})|| \cdot ||{\bf z'}|| \nonumber\\
 &=&  \lambda_j|| {\bf Q}_n^H {\bf a}(\theta_{Ij})||\cdot  || {\bf w}||\nonumber\\
&\ll&   ||{\bf a}(\theta_{Ij})|| \cdot || {\bf w}||.\nonumber
\end{eqnarray}
This means, the beamformer ${\bf \hat w}$ produced by SSC-DL guarantees a suppression level on the signal component ${\bf a}(\theta_{Ij})$ with a response below $\lambda_j$,
whose value is roughly inversely proportional to the root-mean strength of the interference ${\bf a}(\theta_{Ij})$ present in the autocorrelation matrix ${\bf R}_y$.

On the other hand, a beamformer obtained by the conventional DL method, ${\bf w}_{DL}$, does not have such a guarantee.
Adjusting the diagonal loading factor $\gamma$ in the DL method has the effect of reducing attenuation of any signal components from directions other than $\theta_0$, and is, therefore, helpful in relieving undesired suppression of the target signal from angle $\theta_d\neq\theta_0$.
However, it also reduces suppression on any interference signal from angle $\theta_{Ij}$.
By merely adjusting $\gamma$, the DL method is unable to distinguish interference signals from the desired signal and apply different suppression levels on them. 
As the DL objective function (\ref{eq:DL}) does not explicitly prohibit the beamformer to contain signal components in ${\cal W}^\perp$, ${\bf w}_{DL}$ may comprise components from both ${\cal W}$ and ${\cal W}^\perp$.
Assume ${\bf w}_{DL}$ is decomposed as
$ {\bf w}_{DL} = {\bf Q}_s {\bf z}_s + {\bf Q}_n {\bf z}_n$. 
Then, almost certainly ${\bf w}_{DL}$ has a nonzero (though small) component in ${\cal W}^\perp$ (i.e., $||{\bf z}_n||>0$).
Then $ | {\bf a}^H(\theta_{Ij}) {\bf \hat w}| < |{\bf a}^H(\theta_{Ij}){\bf w}_{DL}|$ is very likely to happen, as long as $||{\bf z}_n||/ ||{\bf z}_s|| \gg \lambda_j$, which is usually true, especially when ${\bf a}(\theta_{Ij})$ represents a strong interference source. { In Section \ref{secsec:ProjectionRatio}, we will provide some numerical results to demonstrate that beamformers obtained by conventional DL methods have a larger proportion of energy in the interference domain than those obtained by SSC-DL}.
The observation prevents DL methods from having a good overall SINR performance.  

\section{Derivation of DL factor for SSC-DL}
\label{sec:derivation}
In this section, we present the formulation to determine the DL factor $\gamma$ for the SSC-DL method. The derivation basically follows the flow used in \cite{vincent2004}, with a major difference that the matrix ${\bf W}_M$ is involved here due to the subspace constraint. For notational simplicity, we use ${\bf a} = {\bf a}(\theta_{d})$, $\bar{\bf a} = {\bf a }(\theta_{0})$, and $\sigma_{n}^2$ is noise variance. Our goal is to determine $\gamma$ that maximizes the output SINR, i.e., 
\begin{align}\label{maxsinr}
    \gamma_{opt} = \arg\max_{\gamma}\frac{P\left | \hat{{\mathbf w}}_{M}^{H} {\bf a}(\theta_{d}) \right |^2}{\hat{{\mathbf w}}_{M}^{H}{\bf R}_{I+n}\hat{{\mathbf w}}_{M}}.
\end{align}
Substituting (\ref{ssc-dl}) into the objective function in (\ref{maxsinr}), and with some efforts similar to that in \cite{vincent2004}, it can be shown that:
\begin{align}\label{SINR_newform}
&\frac{P\left | \hat{{\mathbf w}}_{M}^{H} {\bf a}(\theta_{d}) \right |^2}{\hat{{\mathbf w}}_{M}^{H}{\bf R}_{I+n}\hat{{\mathbf w}}_{M}} = \frac{1}{P({\bf a}-\alpha\bar{{\bf a}})^H{\bf W}_{M}{\bf Z}{\bf W}_{M}^H({\bf a}-\alpha\bar{{\bf a}})} 
\end{align}
where the value $\alpha$ and the matrix ${\bf Z}$ are defined as
\begin{equation}\label{alpha}
 \alpha = \frac{1+P{\bf a}^H{\bf W}_{M}({\bf W}_{M}^H({\bf R}_{I+n}+\gamma{\bf I}){\bf W}_{M})^{-1}{\bf W}_{M}^H{\bf a}}{P{\bf a}^H{\bf W}_{M}({\bf W}_{M}^H({\bf R}_{I+n}+\gamma{\bf I}){\bf W}_{M})^{-1}{\bf W}_{M}^H\bar{{\bf a}}}
\end{equation}

\begin{align}\label{Z}
\mbox{and }~~~{\bf Z} &=& ({\bf W}_{M}^H({\bf R}_{I+n}+\gamma{\bf I}){\bf W}_{M})^{-1}{\bf W}_{M}^H{\bf R}_{I+n}{\bf W}_{M}\nonumber\\
&&\cdot({\bf W}_{M}^H({\bf R}_{I+n}+\gamma{\bf I}){\bf W}_{M})^{-1}.
\end{align}
Define the function $f(\gamma)$ as the denominator of the right-hand side of Eq. (\ref{SINR_newform}):
\begin{equation}\label{f}
f(\gamma) = ({\bf a}-\alpha\bar{{\bf a}})^H{\bf W}_{M}{\bf Z}{\bf W}_{M}^H({\bf a}-\alpha\bar{{\bf a}}).
\end{equation}
Then minimizing $f(\gamma)$ is equivalent to the problem in (\ref{maxsinr}).
To simplify the expressions of $\alpha$ and ${\bf Z}$, we assume that ${\bf W}_M$ is ``almost orthogonal'' to any interfering signal ${\bf a}(\theta_j)$ (see Section \ref{subsec:interf_suppr}) and recognize that a common term in $\alpha$ and ${\bf Z}$ can be approximated as
\begin{eqnarray}
&& {\bf W}_M({\bf R}_{I+n}+\gamma{\bf I}){\bf W}_M \nonumber\\
&=& {\bf W}_M\left[\sum_{j=1}^J x_j(t){\bf a}(\theta_j){\bf a}^H(\theta_j)+(\sigma_n^2+\gamma){\bf I}\right]{\bf W}_M \nonumber\\
&\approx& (\sigma_n^2 + \gamma){\bf W}_M^H{\bf W}_M = (\sigma_n^2 + \gamma){\bf V}{\boldsymbol\Sigma}^2_M{\bf V}^H  \label{near_orthogonal}
\end{eqnarray}
where the matrix ${\bf W}_M$ is assumed to have the form of singular value decomposition
$ {\bf W}_M = {\bf Q}_s{\boldsymbol\Sigma}_M{\bf V}^H$.
Note that ${\bf Q}_s$ forms an orthonormal basis of the constraint subspace in SSC-DL.
Then, using (\ref{Z})(\ref{near_orthogonal}), we have
\begin{align}\label{WZW}
{\bf W}_M{\bf Z}{\bf W}_M^H \approx \frac{\sigma_n^2}{(\sigma_n^2 + \gamma)^2}{\bf Q}_s{\bf Q}_s^H.
\end{align}
Similarly, Eq. (\ref{alpha}) can be simplified to
\begin{align}\label{n_alpha}
&\alpha \approx  \underbrace{\frac{\sigma_n^2 + P {\bf a}_n^H {\bf a}_n}{P {\bf a}_n^H \bar{\bf a}_n}}_{\triangleq\alpha_0} + 
{\gamma}\cdot\underbrace{\frac{1}{P {\bf a}_n^H \bar{{\bf a}}_n}}_{\triangleq\alpha'_0} 
\end{align}
where ${\bf a}_{n} = {\bf Q}_s^H{\bf a}$ and $ \bar{{\bf a}}_{n} = {\bf Q}_s^H\bar{{\bf a}} $ represent projections of desired and estimated steering vectors ${\bf a}$ and $\bar{\bf a}$ onto the constraining subspace. Then, using (\ref{WZW}), (\ref{n_alpha}), the function (\ref{f}) can be rewritten as
\begin{align}\label{fgamma}
f(\gamma) = \frac{\sigma_n^2}{(\sigma_n^2 + \gamma)^2} \|{\bf b}+\gamma{{\bf b}}'\|^2,  
\end{align}
where 
\begin{align}\label{algebra}
&{\bf b} = {\bf Q}_s^H ({\bf a}-\alpha_0\bar{{\bf a}}) = {\bf a}_{n} - \alpha_{0}\bar{{\bf a}}_{n} \\
&{\bf{b}}' = -{\bf Q}_s^H{\alpha_0}'\bar{{\bf a}} = -{\alpha_0}'\bar{{\bf a}}_{n} \nonumber
\end{align}
Then, by taking derivative on $f$ in (\ref{fgamma}) with respective to $\gamma$, the optimal solution of $\gamma$ can be found to be
\begin{equation}\label{opt_gamma}
\gamma_{opt} = \frac{\|{\bf b}\|^2 - \sigma_n^2\mathbb{R}\left [ {\bf b}^H{{\bf b}}'\right ]}{\sigma_n^2\|{\bf{b}}'\|^2-\mathbb{R}\left [ {\bf b}^H{\bf{b}}'\right ]}
= - (P\|{\bf a}_n\|^2 + \sigma_n^2).
\end{equation}
Note that the norm of ${\bf a}_n={\bf Q}_s^H{\bf a}$ would be very close to that of ${\bf a}$ as long as the  steering vector of the desired signal has most of its energy in the constraining subspace, leading to
\begin{equation}
\|{\bf a}_{n}\|^2 \cong \|{\bf a}\|^2 = N.\label{ana}
\end{equation}
So (\ref{opt_gamma}) becomes
\begin{equation}\label{hatgamma}
\hat{\gamma} = -(\sigma_n^2 + P\|{\bf a}\|^2) =  -(\sigma_n^2 + P\cdot N).
\end{equation}
This result happens to be identical to that given in \cite{vincent2004}. However, the validity of this optimal solution is subject to two assumptions given in (\ref{near_orthogonal}) and (\ref{ana}). As in SSC-DL, these two assumptions are satisfied, the choice of $\gamma$ presented in (\ref{eq:gamma}) is justified.

\section{Simulation}\label{sec:sim}
In this section, we report the simulation results of the proposed SSC-DL algorithm. 
For comparison purpose, we conducted experiments on several related beamforming algorithms, including (a) MVDR\cite{capon1969}, (b) MVDR with DL \cite{bd1998, li2003, du2010}, (c) a quadratically constrained method with DoA beamwidth widening in \cite{chen2007}, and (d) an autocorrelation matrix reconstruction method in \cite{huang2015}. 
These four methods are labeled as ``MVDR'', ``MVDR-DL'', ``[Chen2007],'' and ``[Huang2015],'' respectively, in the following discussions. 
The optimal result of MVDR method with perfect $\mathbf{R}_{y}$ and without any DoA error is also reported and denoted as ``Optimal,'' which serves as the theoretical limit for all methods.

\subsection{Simulation Setup}
\label{subsec:simsetup}
{In our simulation}, the uniform linear array is employed with $N=10$ omni-directional sensors whose inter-element distance is half the wavelength: $\lambda/2$. 
Except for ``Optimal'', all methods use $K = 100$ snapshots to construct the $\hat{\bf{R}}_y$ matrix defined in (\ref{eq:hatRy}). 
The actual DoA of the desired signal is set to $\theta_d=0^{\circ}$ and the estimated DoA is $\theta_0=2.5^{\circ}$. 
Two interference sources are radiating through the angles $\theta_{I1} = -20^\circ$ and $\theta_{I2} = 30^\circ$ (i.e., $J=2$). Then, ($\ref{eq:yt}$) is rewritten as
\begin{equation*}\label{eq20}
\mathbf{y}(t)=x_{0}(t)\mathbf{a}(\theta_{d})+x_{1}(t)\mathbf{a}(\theta_{I1})+x_{2}(t)\mathbf{a}(\theta_{I2})+\mathbf{n}(t).
\end{equation*}
{The $x_{0}(t)$, $x_{1}(t)$, $x_{2}(t)$, and $\mathbf{n}(t)$ are zero mean wide-sense independent stationary random processes and satisfy}: 
\begin{eqnarray*}
E[|x_{0}(t)|^{2}]=\sigma_{x}^{2}= P, && 
E[|x_{1}(t)|^{2}]=\sigma_{I_{1}}^{2}=10^4(40\mbox{dB}) \\ 
E[\mathbf{n}(t)\mathbf{n}^{H}(t)]=\sigma_{n}^{2}\mathbf{I}_{N}, &&
E[|x_{2}(t)|^{2}]=\sigma_{I_{2}}^{2}=10^2(20\mbox{dB})
\end{eqnarray*}
where $\sigma_{n}^{2}$ = 1. The perfect covariance matrix is expressed as 
\begin{equation}
\begin{aligned}\label{eq21}
\mathbf{R}_{y}&=P\mathbf{a}(\theta_{d})\mathbf{a}^{H}(\theta_{d})+\sigma_{I_{1}}^{2}\mathbf{a}(\theta_{I1})\mathbf{a}^{H}(\theta_{I1})\\
&+\sigma_{I_{2}}^{2}\mathbf{a}(\theta_{I2})\mathbf{a}^{H}(\theta_{I2})+\mathbf{I}_{N}.
\end{aligned}
\end{equation}
All beamformers know that the estimated DoA $\theta_0$ has an error up to $\pm 4^\circ$ (i.e., $\theta_1 = -1.5^\circ$ and $\theta_2 = 6.5^\circ$) \cite{chen2007}. 
We choose $M=5$ for the proposed SSC-DL method. 
For both the MVDR-DL and the proposed SSC-DL methods, the DL factor is chosen as $\gamma = -(\sigma_n^2 + P \cdot N)$ \cite{vincent2004}. 
For [Chen2007], we follow their Algorithm 2 with the initial DL factor set as $\gamma=1$, and with step angles $\zeta_1=1^\circ, \zeta_2 = 2.5^\circ, $ and $ \zeta_3 = 4.5^\circ$ \cite{chen2007}. 
For [Huang2015], we follow the implementation steps reported in \cite{huang2015} by the estimated DoA of interference signals are $\bar{\theta}_{I1} = -17.5^{\circ}$ and $\bar{\theta}_{I2} = 32.5^{\circ}$, and thus $\Theta_{int}$ = $[\bar{\theta}_{I1}-4,\bar{\theta}_{I1}+4]$ $\cup$ $[\bar{\theta}_{I2}-4,\bar{\theta}_{I2}+4]$; the spherical uncertainty constant $\epsilon = 0$, the number of sampling points in $\Theta_{int}$ is $Q$ = 40, and the number of sampling points on the surface of the spherical uncertainty set is $L$=0 \cite{huang2015}. 
The results of above beamforming methods, to be reported below, are the average of 200 Monte Carlo simulations with $x_{0}(t)$, $x_{1}(t)$ and $x_{2}(t)$ having the standard normal distribution.

\subsection{Magnitude Response Versus Incoming Angle}
\label{subsec:caponspectrum}
\begin{figure}[htp]
\centering \centerline{
\includegraphics[width=0.5\textwidth,height=7cm]{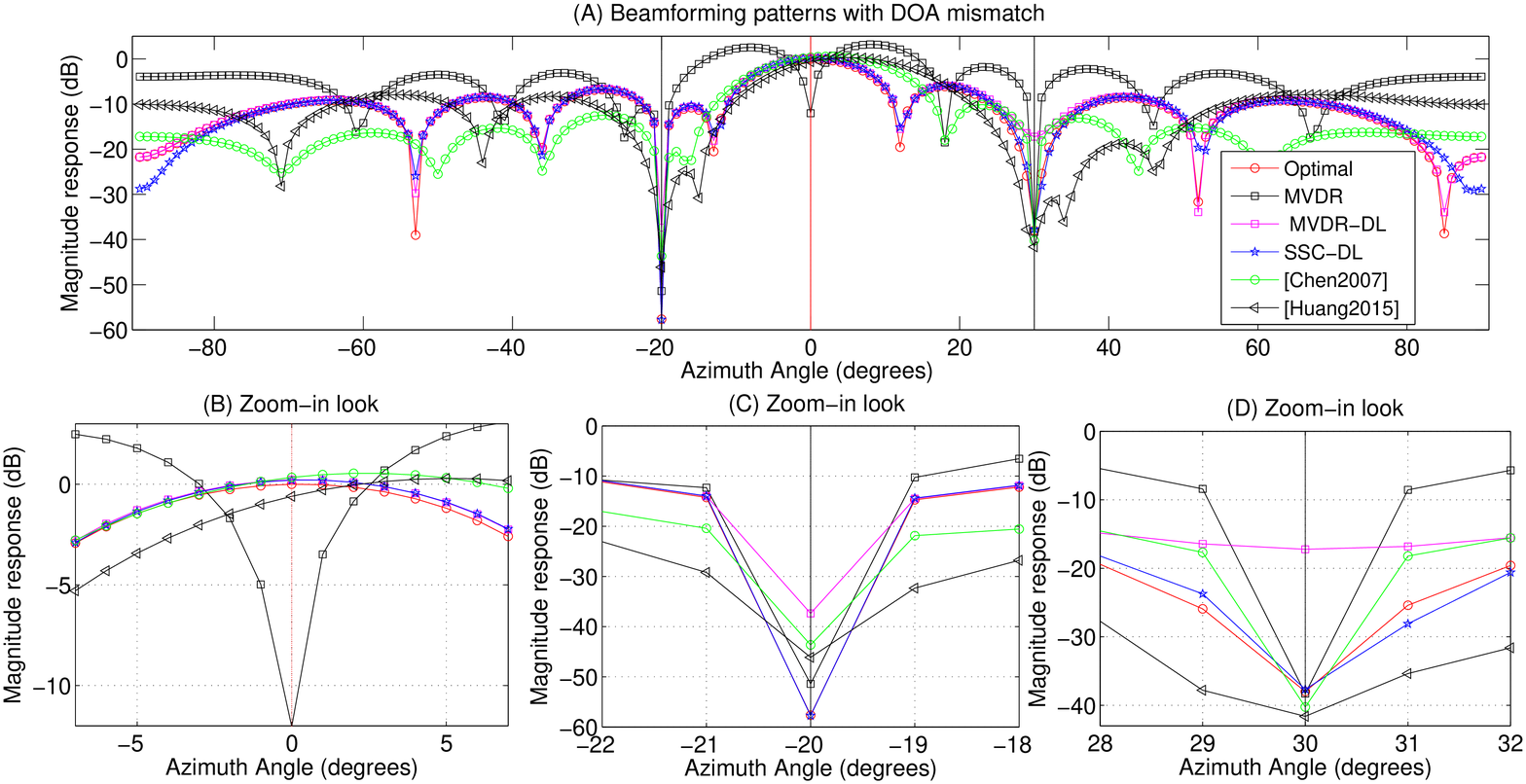}}
\caption{Magnitude response $|\mathbf{w}^H\mathbf{a}(\theta)| $ versus angle $\theta$}
\label{fig4}
\end{figure}

Figure \ref{fig4}(A) shows the beamforming patterns of all methods by plotting the magnitude response versus azimuth angles in degrees. 
The SNR $P/\sigma_n^2$ is set to 10 dB in this plot.
Since the desired signal DoA is at $\theta_d=0^{\circ}$ (red vertical line), and the two interference sources are at $\theta_{I1} = -20^\circ$ and $\theta_{I2} = 30^\circ$(black vertical lines), a successful beamforming technique should give high magnitude response around $0^{\circ}$ and low magnitude response around $-20^\circ$ and $30^\circ$. 
Figures \ref{fig4}(B)(C)(D) display zoom-in looks of Figure \ref{fig4}(A) at angles around $\theta_d$, $\theta_{I1}$, and $\theta_{I2}$, respectively.
From Figure \ref{fig4}(B), we note that due to the DoA mismatch, MVDR gives low magnitude response at $\theta_d=0^{\circ}$ while all  other methods are robust to such a mismatch and generate close-to-unity magnitude responses at the desired DoA. 
From Figures \ref{fig4}(C)(D), SSC-DL shows effective interference signal suppression performance, especially at the angle of $-20^\circ$ where a strong interference resides. 
At the angle of $30^\circ$, SSC-DL does not have a lower response than [Chen2007] and [Huang2015].
But the interferer at $30^\circ$ is 20dB weaker than that at $-20^\circ$, so SSC-DL turns out to still have the largest SINR among all methods, as will be validated later in Section \ref{sub:SINRvsSNR}.
Another interesting observation is the interference signal magnitude responses  of SSC-DL are very close to that of Optimal. 

Further, we compare the main lobes of beamformers produced by SSC-DL and [Chen2007].
Both methods are given the information that the desired angle $\theta_d$ is within the range $[\theta_1, \theta_2] = [-1.5^\circ, 6.5^\circ]$.
The method [Chen2007] uses a strategy to force magnitude responses for all angles in $[-1.5^\circ, 6.5^\circ]$ to be greater than or equal to unity, as can be seen from Figure \ref{fig4}(B).
The SSC-DL, on the other hand, has a narrower beamwidth, with large magnitude responses only at around the desired signal's true DoA.
This suggests setting explicit constraints on maintaining magnitude responses at the range of all possible angles is not necessary.

\subsection{Output SINR versus SNR}

\label{sub:SINRvsSNR}
\begin{figure}[htp]
\centering \centerline{
\includegraphics[width=0.5\textwidth]{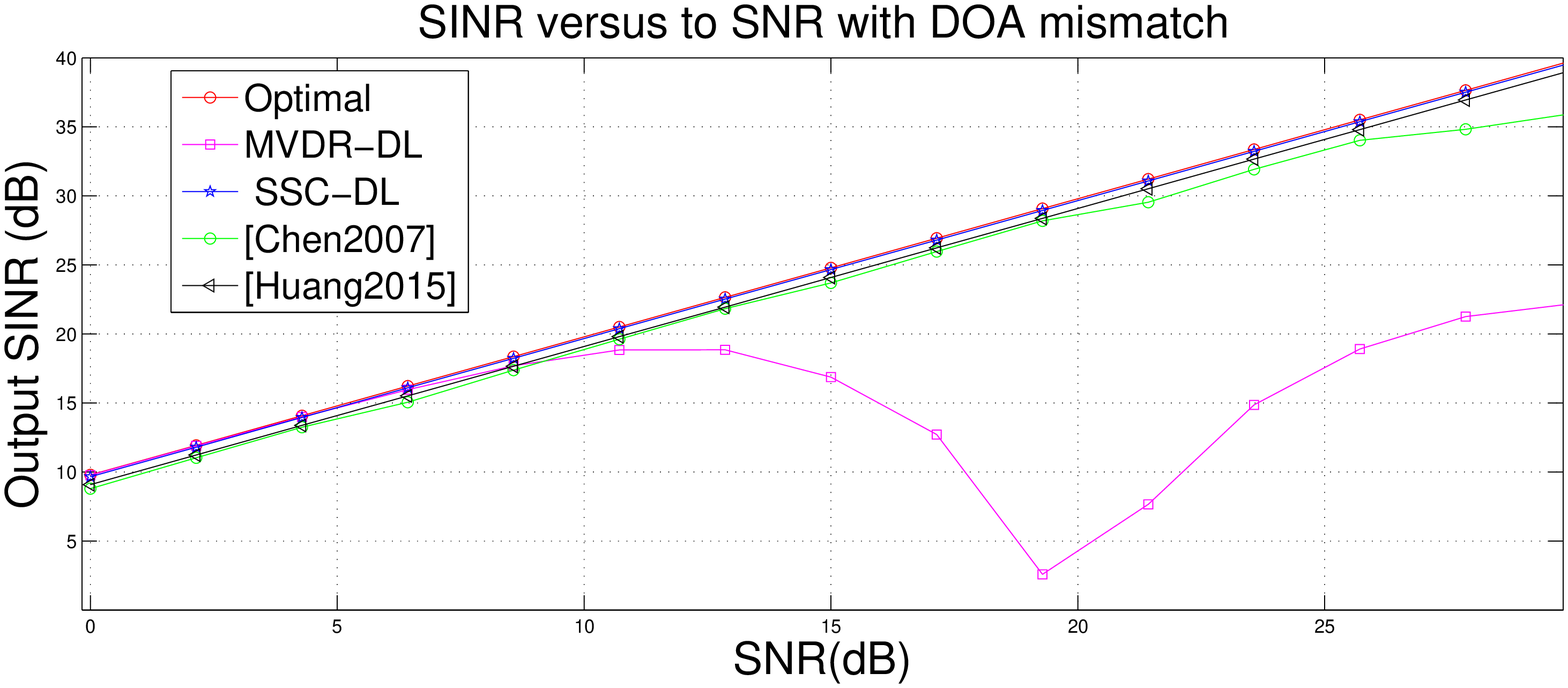}}
\caption{Output SINR versus SNR  }
\label{fig10}
\end{figure}

Figure \ref{fig10} shows the correlation of input SNR ($P/\sigma_n^2$) to the output SINR. 
The plot is generated by fixing $\sigma_n^2, \sigma_{I1}^2$, and $\sigma_{I2}^2$ and by adjusting $P$ to reflect different input SNR levels.
For this set of simulation, given an input SNR, a higher output SINR indicates a better performance. 
From Figure \ref{fig10}, we first observe that MVDR-DL performs poorly when the power of the desired signal is close to the power of some interference (e.g., SNR$\approx 20$ dB); this is caused by the fact that the interference signal is not effectively suppressed, and thus the output SINR is degraded. 
It is clear that SSC-DL successfully overcomes such a limitation by using the subspace technique. 
Moreover, we note that SSC-DL not only outperforms all other algorithms across different input SNR levels between 0 to 30 dB but also achieves a performance very close to the Optimal benchmark.

\subsection{Output SINR versus snapshot number}
\label{subsec:SINRvsK}
\begin{figure}[htp]
\centering \centerline{
\includegraphics[width=0.5\textwidth]{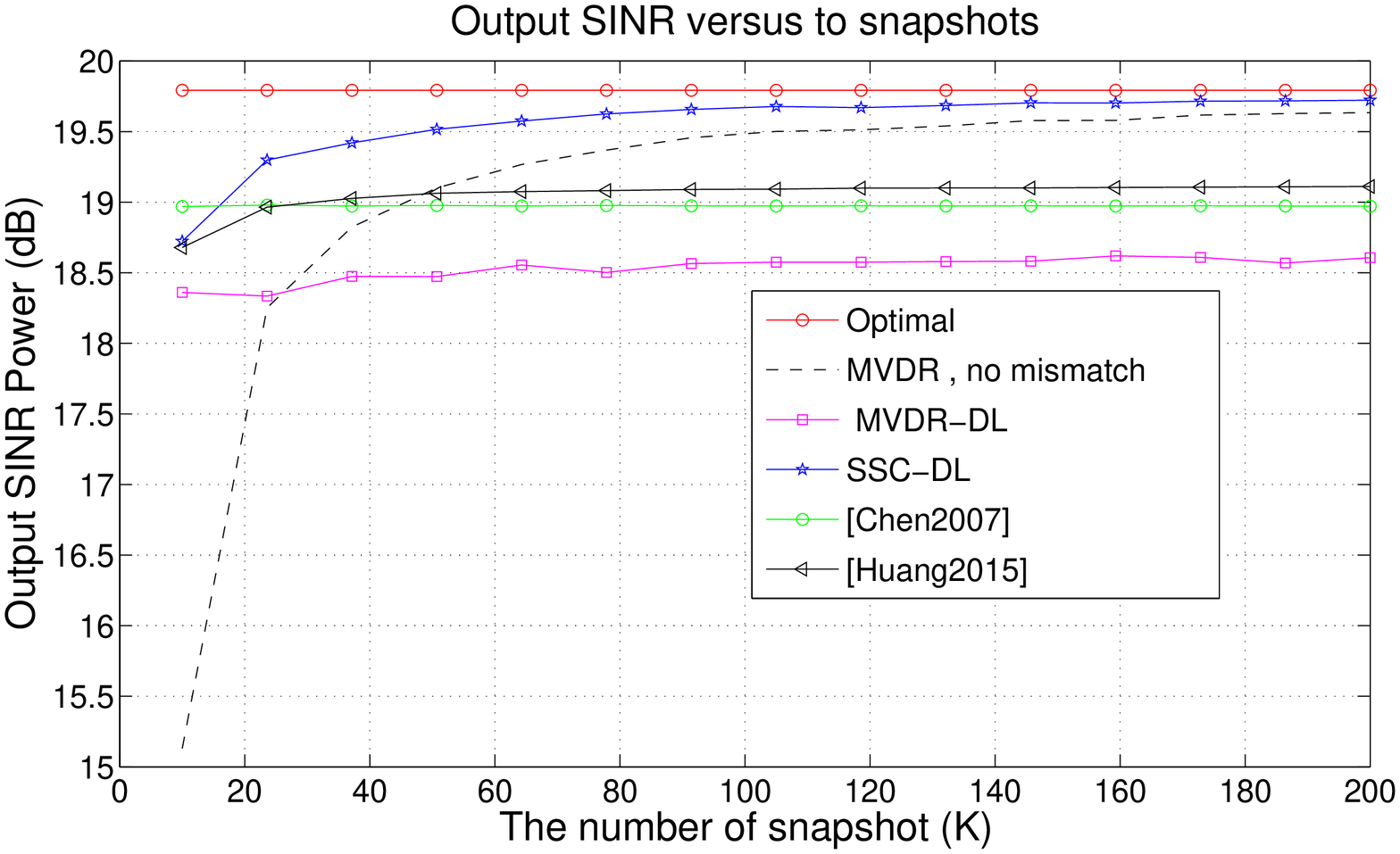}}
\caption{Output SINR versus snapshot number ($K$)}
\label{fig:SINRvsK}
\end{figure}
In Figure \ref{fig:SINRvsK}, we investigate the correlation of the number of snapshots $K$ used in producing $\hat{\bf R}_y$ defined in (\ref{eq:hatRy}) and the output SINR performance of all beamforming methods. 
In addition, a new curve labeled ``MVDR, no mismatch'' is added to show the performance of method \cite{capon1969} without DoA mismatch (i.e., resetting $\theta_0=0^\circ$). 
The SNR is set to 10dB in this plot. 
Inheriting the good property of MVDR-DL against the finite sample effect, SSC-DL also has a good performance with a limited number of snapshots. 
From Figure \ref{fig:SINRvsK}, [Chen2007] has the best performance when very few snapshots are available (e.g., $K<10$). 
When $K >  20$, SSC-DL achieves the best performance among all beamforming methods. 
The output SINR of SSC-DL appears to be very close to the performance limit depicted in ``Optimal'' when $K>160$, but it is not supposed to approach the limit as $K \to \infty$. 
Instead, the curve ``MVDR, no mismatch'' shall approach the limit, as the estimated autocorrelation matrix $\hat{\bf R}_y$ approaches ${\bf R}_y$ as $K\to\infty$. 
Therefore, it would eventually exceed the curve ``SSC-DL.''
Although not shown explicitly in Figure \ref{fig:SINRvsK}, our numerical results indicate that this happens only when $K \geq 12,000$. 
This shows that the proposed SSC-DL not only is robust to DoA mismatch but also is robust to a limited number of available snapshots.

\subsection{The Projection Ratio of Interference signals}
\label{secsec:ProjectionRatio}
\begin{figure}[h]
\centering \centerline{
\includegraphics[width=0.5\textwidth]{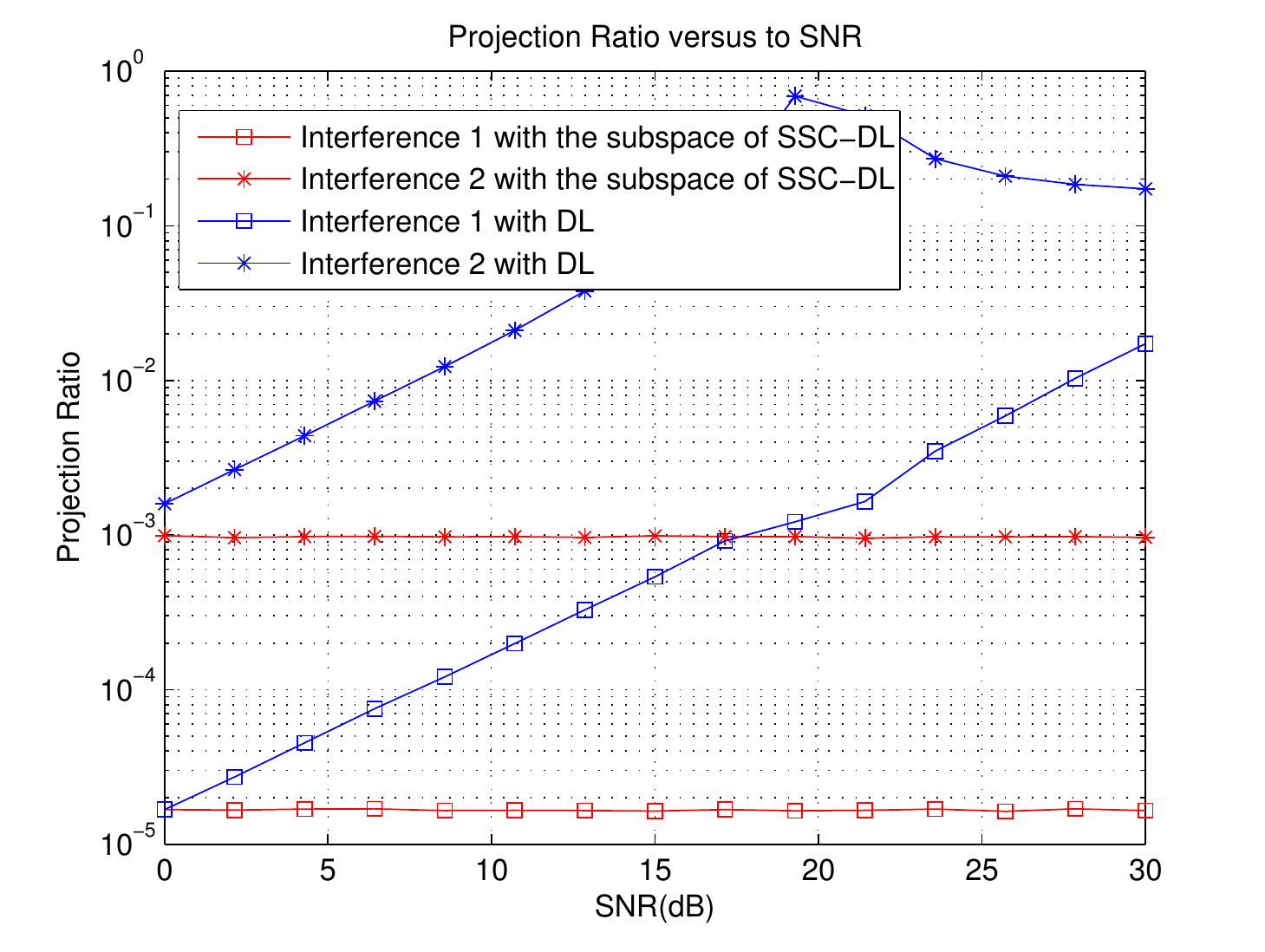}}
\caption{The Projection Ratio of Interference signals versus SNR  }
\label{fig:prvssnr}
\end{figure}

\begin{figure}[h]
\centering \centerline{
\includegraphics[width=0.5\textwidth]{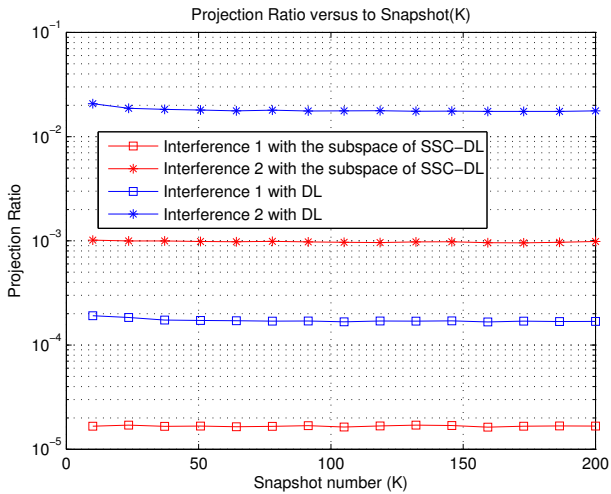}}
\caption{The Projection Ratio of Interference signals versus snapshot number ($K$).  }
\label{fig:prvssnapshot}
\end{figure}

{ In this section, we present the simulation results of projection ratio to validate that the estimated weight vector using the constructed subspace in SSC-DL can more effectively suppress interference than the one estimated by the conventional DL method. The projection ratio for SSC-DL is calculated by the ratio of projections of the steering vector of the interference signal on ${\bf Q}_{s} $ and ${\bf Q}_{n}$:
\begin{align}\label{eq:sscdlpr}
\mbox {Projection Ratio for SSC-DL} =  \frac{\left\| {\bf Q}_{s}{\bf Q}^H_{s}{\bf a}(\theta_{Ij})\right\|}{\left\| {\bf Q}_{n}{\bf Q}^H_{n}{\bf a}(\theta_{Ij})\right\|},
\end{align}
where ${\bf Q}_s$ and ${\bf Q}_n$, as defined in Section \ref{subsec:interf_suppr}, are matrices whose columns form orthonormal bases for the subspaces ${\cal W}$ in SSC-DL and ${\cal W}^\perp$, respectively.} Notably, the projection ratio for SSC-DL indicates the ratio of the closeness of ${\bf a}(\theta_{Ij})$ to ${\bf Q}_s$ over that of ${\bf a}(\theta_{Ij})$ to ${\bf Q}_n$. A smaller projection ratio indicates that ${\bf a}(\theta_{Ij})$ is closer to the matrix ${\bf Q}_n$, which also suggests that ${\bf Q}_s$ of SSC-DL can suppress the interference more effectively.  
{Next, to compute the projection ratio for DL, we project the steering vector of interference signal on ${\bf w}_{DL}$ (the weight vector estimated by the DL method) and a matrix ${\bf Q}_{l}$ whose columns form an orthonormal basis for the null space of ${\bf w}_{DL}$. Then the projection ratio for DL is defined by:
\begin{align}\label{eq:dlpr}
\mbox {Projection Ratio for DL} =  \frac{\left\|  \frac{{\bf a}^H(\theta_{Ij}){\bf w}_{DL}}{\left \| {\bf w}_{DL} \right \|^2} {\bf w}_{DL} \right\|}{\left\| {\bf Q}_{l}{\bf Q}^H_{l}{\bf a}(\theta_{Ij})\right\|}.
\end{align}
Similar to the projection ratio for SSC-DL in (\ref{eq:sscdlpr}), the projection ratio for DL indicates the ratio of the closeness of ${\bf a}(\theta_{Ij})$ to ${\bf w}_{DL}$ over that of ${\bf a}(\theta_{Ij})$ to ${\bf Q}_l$. A smaller projection ratio indicates that ${\bf a}(\theta_{Ij})$ is closer to the matrix ${\bf Q}_l$, which also suggests that ${\bf w}_{DL}$ can suppress the interference more effectively.

Figure \ref{fig:prvssnr} and Figure \ref{fig:prvssnapshot} show the results of projection ratios versus SNR(dB) and snapshot number ($K$), respectively. Note that Figure \ref{fig:prvssnapshot} is showed with fixed SNR = 10 dB. The curves of ``Interference $j$-th with the subspace of SSC-DL'' and ``Interference $j$-th with DL'' are computed by (\ref{eq:sscdlpr}) and (\ref{eq:dlpr}), respectively. From these two figures, the curve ``Interference 1 with the subspace of SSC-DL'' yields smaller projection ratios than the curve ``Interference 1 with DL'', confirming that the proposed SSC-DL can estimate the weight vectors in a more effective way to suppress interference 1 than the DL method. Similarly, we note that SSC-DL also provides more suppression for interference 2, again confirming that the proposed SSC-DL outperforms the DL method in terms of suppressing interference signals to achieve better output SINR performance. }
\label{pr}

\subsection{SINR versus subspace dimension}
\label{subsec:SINRvsM}
\begin{figure}[h]
\centering \centerline{
\includegraphics[width=0.5\textwidth]{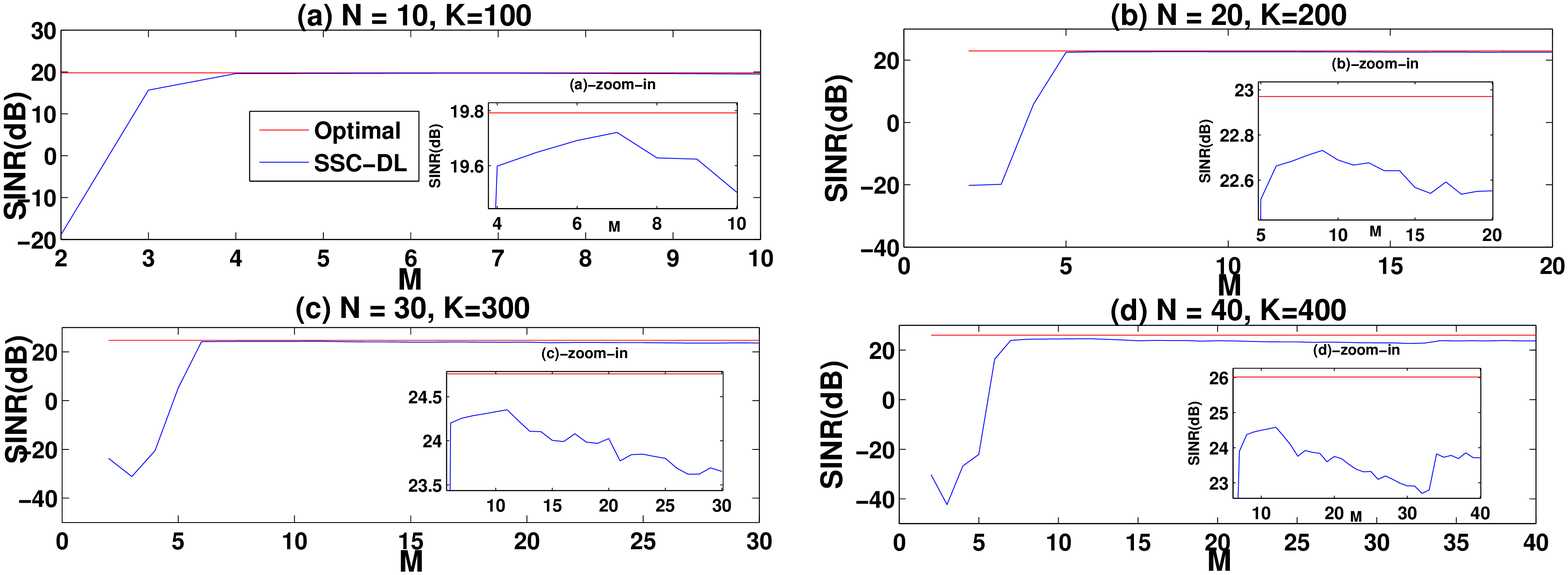}}
\caption{Output SINR versus the subspace dimension $M$.  }
\label{fig:sinrvsm}
\end{figure}

In Figure \ref{fig:sinrvsm}, we study the SINR performance versus the given subspace dimension $M$. As described in Section \ref{subsec:M}, we have learned that $M$ must be large enough to support any possible desired DoA within the range $[\theta_1, \theta_2]$. However, $M$ cannot be too large either: when $M=N$, the SSC-DL degenerates to the conventional DL and has an SINR performance as poor as conventional DL.

\subsection{SINR performance with various DoAs}
\label{subsec:SINRwithVariousDoA}
In this subsection, we display numerical results that show the effectiveness of the proposed SSC-DL method in the environments of various DoAs of the desired signal. 

\begin{figure}[htp]
\centering \centerline{
\includegraphics[width=0.5\textwidth]{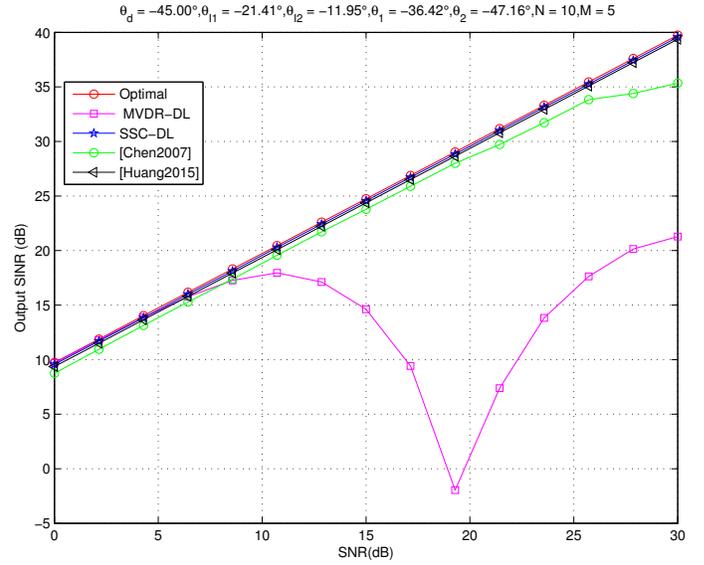}}
\caption{Output SINR versus SNR  when $\theta_d = -45^\circ$. }
\label{fig:snr7}
\end{figure}
\begin{figure}[htp]
\centering \centerline{
\includegraphics[width=0.5\textwidth]{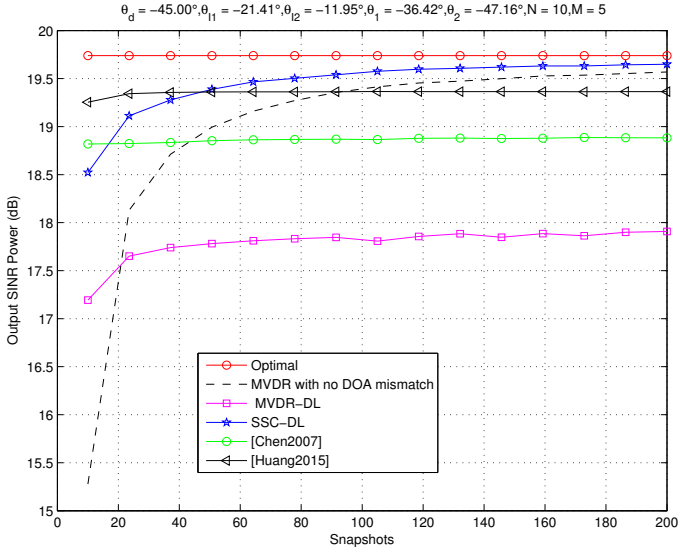}}
\caption{Output SINR versus snapshot number  when $\theta_d = -45^\circ$. }
\label{fig:snapshot7}
\end{figure}

Figures \ref{fig:snr7} and \ref{fig:snapshot7} show the SINR performance of all methods when the DoA of the desired signal is set to $\theta_d = -45^\circ$, and the DoAs of two interference signals are set to $\theta_{I1} = -21.41^\circ$ and $\theta_{I2} = -11.95^\circ$. The estimated DoA of the desired signal is $\theta_0 = \sin^{-1}\left(\sin\theta_d+\sin 2.5^\circ\right)$ and the values of the boundary angles $\theta_1$ and $\theta_2$ are chosen as $\theta_1 = -36.42^\circ$ and $\theta_2 = -47.16^\circ$, following the formulas
$\theta_1 = \sin^{-1}\left(\sin\theta_0-\sin 4^\circ\right) $ and
$\theta_2 = \sin^{-1}\left(\sin\theta_0+\sin 4^\circ\right).$ 
The power levels of the interference are chosen the same as described in Section \ref{subsec:simsetup}:
$E[|x_{1}(t)|^{2}]=\sigma_{I_{1}}^{2}=10^4(40\mbox{dB})$ and 
$E[|x_{2}(t)|^{2}]=\sigma_{I_{2}}^{2}=10^2(20\mbox{dB})$.

In these plots, we find that the proposed SSC-DL remains its superiority in SINR performance to all other methods.
When the number of snapshots is small ($K<10$), the SINR of SSC-DL is slightly lower than those of [Chen2007] and [Huang2015]. 
However, when $K>60$, it is observed that the SSC-DL has the best SINR performance among all methods.

\begin{figure}[htp]
\centering \centerline{
\includegraphics[width=0.5\textwidth]{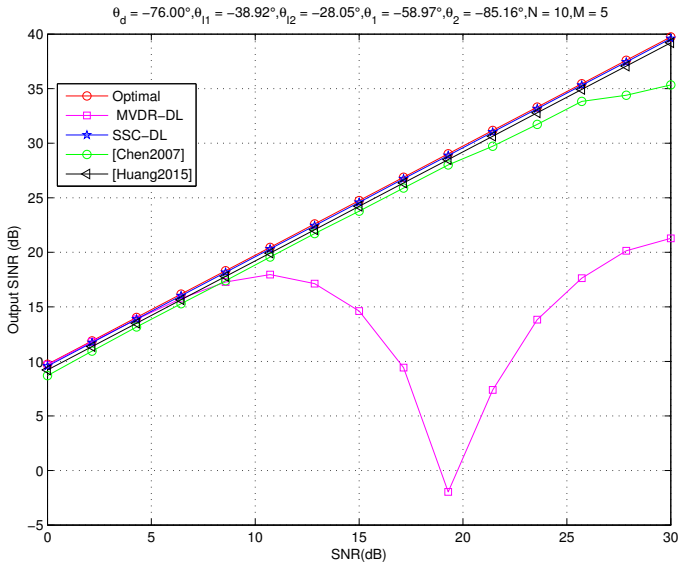}}
\caption{Output SINR versus SNR when $\theta_d = -76^\circ$.}
\label{fig:snr8}
\end{figure}

\begin{figure}[htp]
\centering \centerline{
\includegraphics[width=0.5\textwidth]{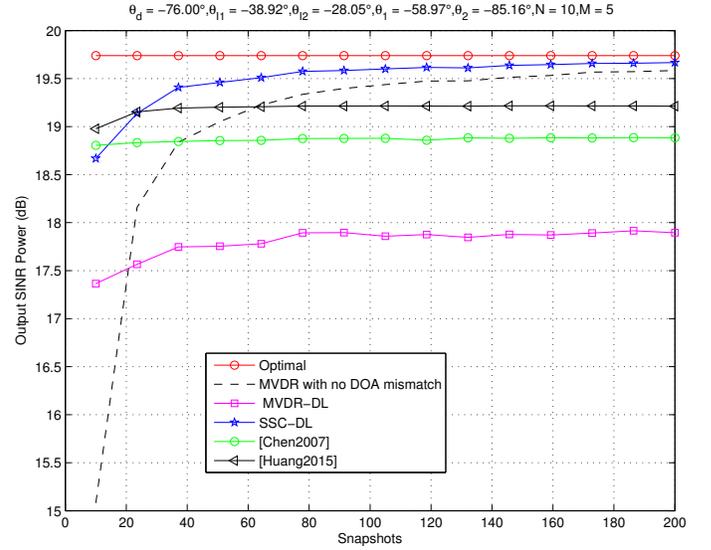}}
\caption{Output SINR versus snapshot number when $\theta_d = -76^\circ$.}
\label{fig:snapshot8}
\end{figure}
Figures \ref{fig:snr8} and \ref{fig:snapshot8} display the SINR performance when the desired DoA is $\theta_d = -76^\circ$. This is an angle that a beamformer based on phase arrays would have a poor resolution. 
The interference signals have DoAs at $\theta_{I1} = -38.92^\circ$ and $\theta_{I2} = -28.05^\circ$.
The estimated DoA of the desired signal is $\theta_0 = -67.92^\circ$
The range of possible desired signal DoA $[\theta_1, \theta_2]$ where $\theta_1 = -58.97^\circ$ and $\theta_2 = -85.16^\circ$, following the formulas
$\theta_1 = \sin^{-1}\left(\sin\theta_0-\sin 4^\circ\right) $ and
$\theta_2 = \sin^{-1}\left(\sin\theta_0+\sin 4^\circ\right).$ 
The power levels of the interference are chosen the same as described in Section \ref{subsec:simsetup}:
$E[|x_{1}(t)|^{2}]=\sigma_{I_{1}}^{2}=10^4(40\mbox{dB})$ and 
$E[|x_{2}(t)|^{2}]=\sigma_{I_{2}}^{2}=10^2(20\mbox{dB})$.

In the simulation plots, we observe that the proposed method still has the best performance among all methods.
The SINR performance of SSC-DL is very close to that of Optimal as long as a sufficient number of snapshots are available ($K>=200$).
When the number of snapshots is small ($K<10$), the proposed method performs slightly worse than [Huang2015], but when $K>30$, it outperforms all other methods.
%These results suggest that the superiority of SSC-DL does not vary.
{Also, according to simulation results with various settings (not present here due to page limit), we found that generally our method would outperform all other existing methods when $K$ is chosen as $\geq 5 \cdot N$ (where $N$ is the number of antennas); when $K$ is chosen as $\geq 20 \cdot N$, the SINR performance roughly reaches its limit.}

\subsection{SINR performance in a larger signal dimension}

\begin{figure}[htp]
\centering \centerline{
\includegraphics[width=0.5\textwidth]{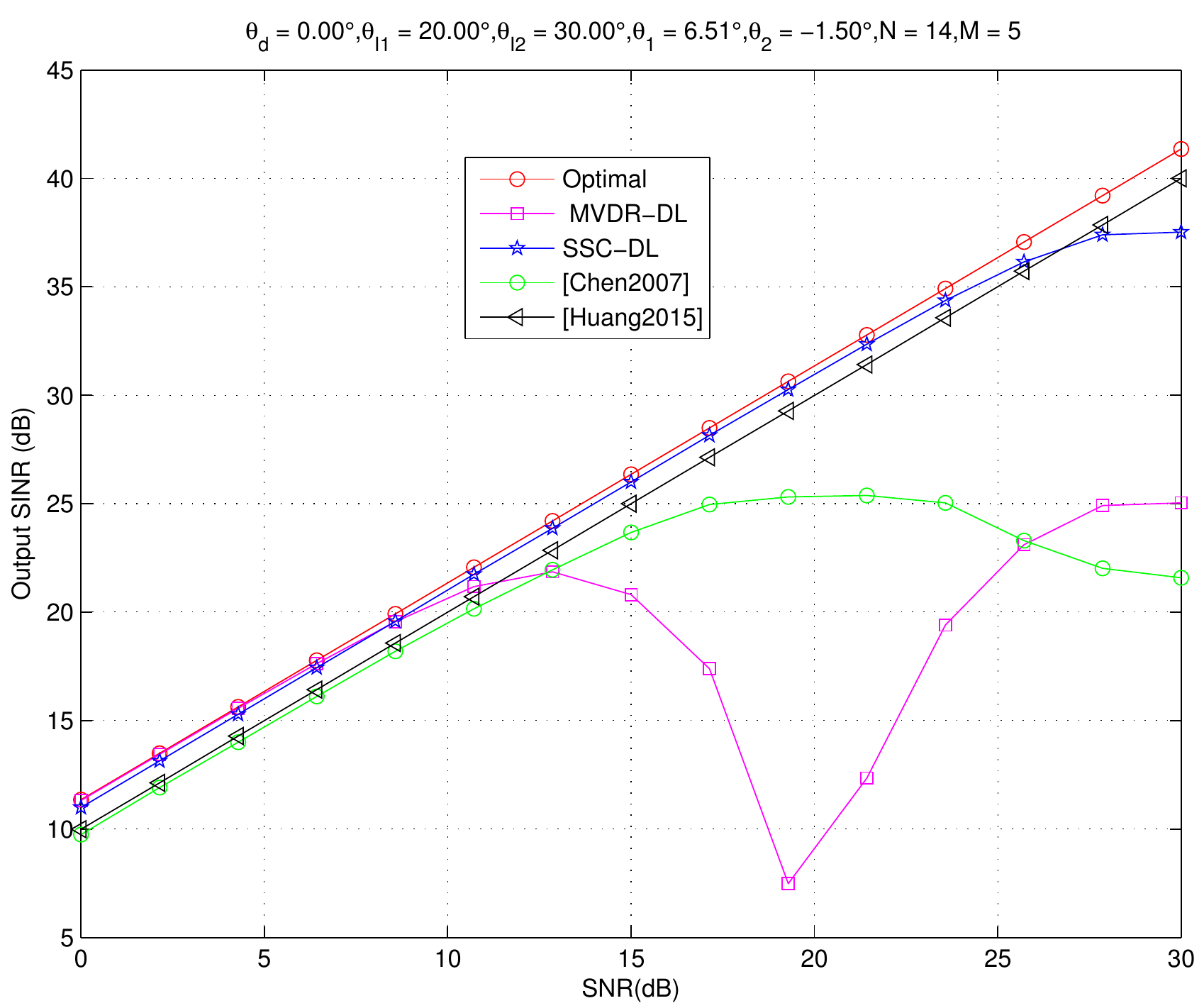}}
\caption{Output SINR versus SNR when $N=14$. }
\label{fig:snr14}
\end{figure}

\begin{figure}[htp]
\centering \centerline{
\includegraphics[width=0.5\textwidth]{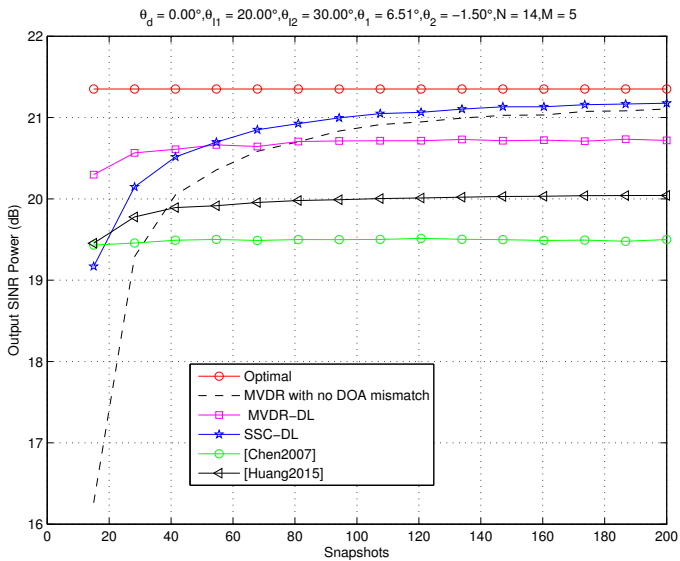}}
\caption{Output SINR versus Snapshot when $N=14$. }
\label{fig:snapshot14}
\end{figure}

Figures \ref{fig:snr14} and \ref{fig:snapshot14} show the SINR performance when the number of antennas increases to $N = 14$.
The desired DoA is $\theta_d = 0^\circ$. Other simulation parameters are the same as those in Section \ref{subsec:simsetup}. 
In this case, the SINR performance of [Chen2007] degrades in the high SNR region since the method \cite{chen2007} requires the parameter $N$ to be upper bounded by a factor in order to have a satisfactory performance.
The proposed SSC-DL still has an excellent performance and outperforms all other methods.

\section{Conclusion} \label{sec:conclusion}
In this paper, we proposed a novel yet simple method, namely subspace-constrained diagonal loading (SSC-DL), to deal with the DoA mismatch issue in a beamforming system. 
The proposed method imposes an additional subspace constraint on the conventional diagonal loading (DL) method to provide better interference suppression and has a closed-form solution. 
The proposed method maintains a high magnitude response to the desired DoA but does not have to widen the beamwidth like other existing approaches.
The signal-to-interference-plus-noise ratio (SINR) performance, as shown in numerical results, is even approaching the theoretical limit suggested in \cite{capon1969}. 
In the future, an analytical study would be helpful in further understanding the superiority of the proposed method. 
Further, automatic selection of parameters $M$ and $\gamma$ is challenging yet desirable.
In addition, {making the DoA of the desired signal modeled as a random variable} and extending the proposed SSC-DL method to the case of general steering vector errors is also of interest. 

\section*{Statement on Conflict of Interests}
The authors declare that there is no conflict of interests regarding the publication of this paper.

\ifCLASSOPTIONcaptionsoff
  \newpage
\fi
\printbibliography

\end{document}